\DeclareFontFamily{U}{wncy}{}
\DeclareFontShape{U}{wncy}{m}{n}{%
   <5>wncyr5%
   <6>wncyr6%
   <7>wnyr7%
   <8>wncyr8%
   <9>wncyr9%
   <10>wncyr10%
   <11>wncyr10%
   <12>wncyr6%
   <14>wncyr7%
   <17>wncyr8%
   <20>wncyr10%
   <25>wncyr10}{}
\DeclareMathAlphabet{\cyrille}{U}{wncy}{m}{n}
\documentclass[a4paper,11pt]{article}

\usepackage{latexsym}
\usepackage{amsmath,amsfonts,amssymb,amsthm}
\newtheorem{thm}{Theorem}[section]
\newtheorem{lem}[thm]{Lemma}

\newtheorem{cor}[thm]{Corollary}
\newtheorem{prop}[thm]{Proposition}

\newtheorem{rem}[thm]{Remark}

\input epsf

\title{Pr\"ufer codes for hypertrees}
\author{Roland Bacher}

\begin{document}
\maketitle

{\sl Abstract\footnote{
Math. class: Primary: 05C30, Secondary: 05A05, 05A15, 05A19,
05C05, 05C07, 05C65}: We describe a new Pr\"ufer code which works
also for infinite hypertrees. 
}

%fichier prufer1.tex dans polyuniv

\section{Main results}

A famous theorem attributed to Cayley states that there are $n^{n-2}$ finite trees 
with vertices $\{1,\dots,n\}$. 
Pr\"ufer gave in \cite{P} a beautiful proof by constructing a 
one-to-one correspondence between such trees and elements in the 
set $\mathcal \{1,\dots,n\}^{n-2}$ of all $n^{n-2}$
words of length $n-2$ in the 
alphabet $\{1,\dots,n\}$. More precisely, we obtain 
the Pr\"ufer code of a tree with $n\geq 2$ 
vertices $\{1,\dots,n\}$ by successively pruning
smallest leaves and writing down their neighbours until 
reaching a tree reduced to a unique edge. Selivanov in \cite{Se}
generalized Pr\"ufer's construction to hypertrees.

Pr\"ufer's construction and its subsequent generalizations can be 
succinctly described as \lq\lq pruning trees''. The aim of this
paper is to describe a Pr\"ufer code 
based on a different kind of simplification, star-reduction,
which merges hyperedges until reaching the trivial hypertree consisting of 
a unique hyperedge containing all vertices. The resulting Pr\"ufer
code respects degrees (a vertex of degree
$a$ occurs with multiplicity $a-1$ in the corresponding Pr\"ufer 
word) a property shared with Pr\"ufer's original
construction. Its definition is perhaps slightly less straightforward but
it has the additional feature that it works for infinite trees
and hypertrees, a property which fails to hold in the general case
for the classical construction of Pr\"ufer.

The rest of this paper is organized as follows: Section \ref{sectlabhyp}
recalls briefly the definition of hypertrees. Section \ref{sectPpartition}
describes the Pr\"ufer partition. The Pr\"ufer partition of an ordinary 
rooted tree is trivial and carries no information. It is however a 
necessary ingredient for the Pr\"ufer code of a hypertree with hyperedges
of size larger than $2$. Section \ref{sectgluemap} contains an intuitive
definition of Pr\"ufer codes and reviews a few enumerative results.
Sections \ref{sectleaftype} and \ref{sectinvW} recall the definition of the classical Pr\"ufer
code and its generalization to finite rooted hypertrees. Sections 
\ref{sectlabhyp}-\ref{sectinvW} are expository and contain nothing new.

Sections \ref{sectstarred}-\ref{sectinvW*} describe the Pr\"ufer code
$T\longmapsto (\mathcal P,W^*)$ based on star-reductions (mergings of 
hyperedges). As far as I am aware, this construction has not appeared 
elsewhere in the litterature. 

Section \ref{sectinftrees} generalizes the construction of the Pr\"ufer
code based on star-reductions to infinite trees and hypertrees. 

The final Section \ref{sectexpleSn} illustrates the extension to infinite
trees by an interesting example which gives rise to bijections in the set
$S_n$ of all elements in the symmetric group of $\{1,\dots,n\}$.

\section{Hypergraphs and hypertrees}\label{sectlabhyp}

A \emph{hypergraph} is a pair $(\mathcal V,\mathcal E)$
consisting of a set $\mathcal V$ of \emph{vertices} 
and of a set $\mathcal E$ of \emph{hyperedges}
given by subsets of $\mathcal V$ containing at least two elements.
A hypergraph is \emph{finite} if it contains only finitely many vertices
and finitely many hyperedges.

Except in Sections \ref{sectinftrees} and \ref{sectexpleSn} 
we consider henceforth mainly only finite hypergraphs (and hypertrees) 
consisting of a finite number of vertices and hyperedges.
We denote by $\mathrm{size}(e)$ the cardinality of a hyperedge $e$, 
defined as the number of vertices contained in $e$,
and by $\deg(v)$ the degree of a vertex $v$ 
given by the number of hyperedges containing $v$.
A vertex of degree $1$ is a \emph{leaf}. Both numbers $\mathrm{size}(e)$
and $\deg(v)$ can be arbitrary (perhaps infinite) cardinal numbers.
A hypergraph is \emph{locally finite} if it has only edges of finite size and
vertices of finite degree. 

Two distinct vertices $v,w\in\mathcal V$ of a common hyperedge 
are \emph{adjacent} or \emph{neighbours}.
A \emph{path} of \emph{length $l$} joining two vertices 
$v,w\in\mathcal V$ is a sequence
$v=v_0,v_1,\dots,v_l=w$ involving only
consecutively adjacent (and distinct) vertices. 
A hypergraph is \emph{connected} if two vertices 
can always be joined by a path.  
The \emph{distance} between two vertices $v,w$
of a connected hypergraph 
is the length of a shortest path joining $v$ and $w$.
\emph{Geodesics} are paths $\dots,v_i,\dots,v_j,\dots$ with 
$v_i,v_j$ at distance $\vert i-j\vert$ for all indices $i,j$.
We set $d(v,w)=\infty$ if $v$ and $w$ belong to different connected 
components.

We have
\begin{eqnarray}\label{ineqconhyp}
  \sum_{e\in\mathcal E}\mathrm{size}(e)=\sum_{v\in\mathcal V}\deg(v)\geq n+k-1
\end{eqnarray}
for every connected finite hypergraph with $n$ vertices and $k$ hyperedges.
A connected finite hypergraph is a finite \emph{hypertree} if equality holds in 
(\ref{ineqconhyp}). We have thus
\begin{eqnarray}\label{edgedegineq}
\sum_{e\in\mathcal E}\left(\mathrm{size}(e)-1)\right)=n-1
\end{eqnarray}
and
\begin{eqnarray}\label{degreeequality}
\sum_{v\in\mathcal V}\left(\deg(v)-1)\right)=k-1
\end{eqnarray}
for a finite hypertree with $n$ vertices and $k$ hyperedges.
A hypertree is \emph{trivial}  if it is reduced to a unique hyperedge.
Equivalently, a connected graph is a trivial hypertree 
if all its vertices are leaves.

\begin{prop} (i) We have $k<n$ for a finite hypertree with $n$ vertices and
$k$ hyperedges.

\ \ (ii) Every finite hypertree with $n\geq 2$ vertices contains at least two leaves.
\end{prop}

\noindent{\bf Proof} Since every hyperedge $e$ contains at least 
$\mathrm{size}(e)\geq 2$ vertices we have
$n-1=\sum_{e\in\mathcal E} \left(\mathrm{size}(e)-1\right)\geq k$. This
shows $(i)$. Since $n-1>k-1=\sum_{v\in\mathcal V}\left(\deg(v)-1\right)$
there exists at least two vertices contributing nothing to the 
sum $ \sum_{v\in\mathcal V}\left(\deg(v)-1\right)$.\hfill$\Box$

A connected infinite hypergraph is a \emph{hypertree} if every connected
subgraph induced by a finite number of vertices (with hyperedges 
given by intersections containing at least two vertices of
original hyperedges with the finite subset of vertices under consideration)
is a finite hypertree.

\begin{prop} Two vertices $v,w$ at distance $l$ in a hypertree $T$
are joined by a unique shortest path
$v=v_0,v_1,\dots,v_l=w$ defining a unique sequence $e_1,\dots,e_l$ 
of hyperedges such that $\{v_{i-1},v_i\}\subset e_i$.
\end{prop}

\noindent We leave the proof to the reader.\hfill $\Box$

%%%%%%%%%%%%%%%%%%%%%%%%%%%%%%%%%%%%%%%%%%%%%%%%%%%%%%%%%%%%%%%%%%%%%%%%%%%%%%%
\section{The  Pr\"ufer partition of a rooted hypertree}
\label{sectPpartition}
A {\emph{rooted hypertree} has a marked root vertex $r$ among its vertices.
The root vertex $r$ induces a \emph{marked vertex}
$e_*$ closest to the root in every hyperedge $e$ of a rooted hypertree.
Removal of the marked vertex 
$e_*$ from $e$ yields the \emph{reduced hyperedge}
$e'=e\setminus\{e_*\}$ consisting of all \emph{unmarked vertices} of $e$. 
No reduced hyperedge contains the root.
We use the notation $\{v_1,\dots,v_{k-1}\}_w$ for a hyperedge 
$e=e'\cup\{e_*\}$ of size $k$
with marked vertex $e_*=w$ and reduced hyperedge $e'=\{v_1,\dots,v_{k-1}\}$.

\begin{prop}\label{proproot} A non-root vertex $v$ of a rooted hypertree 
is the marked vertex of $\deg(v)-1$ hyperedges. 
The root $r$ is the marked vertex of $\deg(r)$ hyperedges.
\end{prop}

\noindent The proof is left to the reader.
\hfill$\Box$

\begin{cor}\label{corroot} Reduced hyperedges
of a rooted hypertree $T$ with vertices $\mathcal V$, root $r$ 
and $k$ hyperedges
partition the set $\mathcal V\setminus\{r\}$ 
into $k$ non-empty subsets.
\end{cor}

The \emph{Pr\"ufer partition $\mathcal P(T)$} 
of an $r-$rooted hypertree $T$ with vertices $\mathcal V$ is the partition
of $\mathcal V'=\mathcal V\setminus\{r\}$ into 
reduced hyperedges given by Corollary  \ref{corroot}.

The definition of 
Pr\"ufer partitions can easily be generalized to 
arbitrary (not necessarily locally finite) infinite rooted hypertrees.

In the sequel, we will mainly consider hypertrees with non-root vertices
$\mathcal V'$ given by (perhaps infinite) subsets of $\{1,2,\dots\}$. 
The Pr\"ufer partition of such a hypertree is thus either given by a 
partition of $\mathcal V'$ or equivalently by a partition map of 
$\mathcal V'$. 

%%%%%%%%%%%%%%%%%%%%%%%%%%%%%%%%%%%%%%%%%%%%%%%%%%%%%%%%%%%%%%%%%%%%%%%%%%

\subsection{Partition maps}

A map $p:E\longrightarrow E$ of a set $E$ is \emph{idempotent} if
$p=p\circ p$. Equivalently, a map $p:E\longrightarrow E$ is idempotent
if its image $p(E)$ is its set of fix-points.
Idempotent maps of a set $E$ are in one-to-one correspondence with 
partitions of $E$ decorated with a marked element in each part.
Indeed, such a decorated partition gives rise to an idempotent map by
sending each element to the marked element of its part. In the opposite
direction, an idempotent map $p$ gives rise to a partition with marked
elements given by fix-points and parts given by preimages of fix-points.

A set $E$ is well-ordered if $E$ is endowed by an order relation which 
yields a least element in every non-empty subset of $E$. 
A map $p:E\longrightarrow 
E$ of a well-ordered set is \emph{lowering} if $p(x)\leq x$ for all $x$.

Well-ordering a set $E$ selects least elements as the
canonical marked elements in parts of a partition. Partitions of
a well-ordered set $E$ are thus in one-to-one correspondence with 
maps $p:E\longrightarrow E$ which are idempotent and lowering.
We call such a map a \emph{partition map}.

%%%%%%%%%%%%%%%%%%%%%%%%%%%%%%%%%%%%%%%%%%%%%%%%%%%%%%%%%%%%%%%%%%%%%%%%%%%

\subsection{An example of a Pr\"ufer partition
}
\label{sectexpleP}

We consider the finite rooted hypertree $T$ with 
vertices $1,\dots,14$, root-vertex $14$ and $8$ hyperedges given 
by
\begin{eqnarray*}
\{1,10,12\}_8,\{2\}_1,\{3,9\}_4,\{4,7\}_8,
\{5\}_{14},\{6\}_4,\{8,13\}_{14},\{11\}_7
\end{eqnarray*}
where $\{3,9\}_4$ for example represents a hyperedge
of size $3$ with marked vertex $4$
and reduced hyperedge $\{3,9\}$.
\begin{figure}[h]\label{figure1}
\epsfysize=4.5cm
\center{\epsfbox{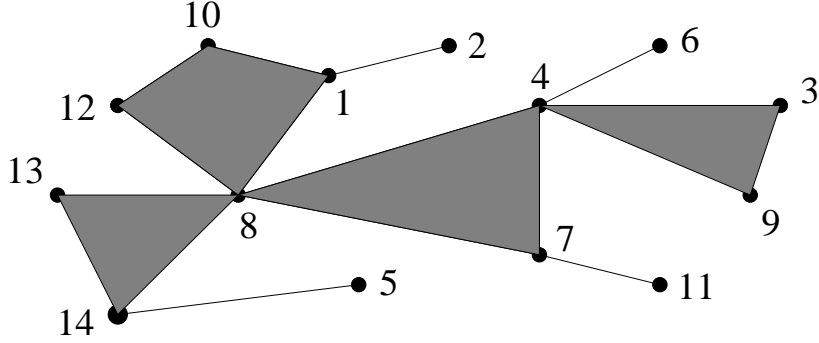}}
\caption{A hypertree}
\end{figure}

Figure 1 shows $T$ with hyperedges represented
by shaded polygons or ordinary edges.

The Pr\"ufer partition $\mathcal P(T)$ of the hypertree $T$ with
root vertex $14$ 
is thus the partition
$$\{1,10,12\}\cup\{2\}\cup\{3,9\}\cup\{4,7\}\cup\{5\}\cup\{6\}\cup
\{8,13\}\cup\{11\}$$
defined by the union of all reduced hyperedges of $T$.
The corresponding partition map is given by
\begin{eqnarray*}
&&p(1,10,12)=1,\ p(2)=2,\ p(3,9)=3,\ p(4,7)=4,\\ 
&&p(5)=5,\ p(6)=6,\ p(8,13)=8,\ p(11)=11\ .
\end{eqnarray*}

\subsection{The spine of a rooted hypertree}
We associate in this short digressional section an ordinary rooted tree
(the spine) to every rooted hypertree.

The \emph{spine} of a rooted hypertree $T$
with root $r$ and vertices $\mathcal V$ 
is the ordinary tree $Sp(T)$ with root $r$, vertices 
$\overline{\mathcal V}=\mathcal P\cup\{r\}$ where elements of $\mathcal P$
are parts involved in the Pr\"ufer partition $\mathcal P$ of $T$
and edges $\{A,B\}$ if an element of $A$ is adjacent to an element of $B$. 
There is an obvious projection
$\pi:\mathcal V\longrightarrow\overline{\mathcal V}$
defined by $\pi(r)=r$ and $\pi(v)=e'$ if $v\in \mathcal V\setminus\{r\}$
is contained in the reduced hyperedge $e'$.
Edges of $Sp(T)$ are in one-to-one correspondence with hyperedges of $T$ 
and are given by $\{e',\pi(e_*)\}$ for a hyperedge $e$ of $T$.

We state the following result without proof:

\begin{prop}\label{propspine} (i) $Sp(T)=T$ if and only if $T$ is an ordinary rooted tree.

\ \ (ii) $d(\pi(v),\pi(w))\leq d(v,w)$ for $v,w\in \mathcal V$ and
$d(\pi(v),\pi(r))=d(v,r)$ for every vertex $v$ in $\mathcal V$.
\end{prop}

\begin{rem} There exist a second natural map which associates a rooted tree 
to every rooted hypertree and which does not modify the set of vertices:
replace every hyperedge $e$ of $T$ with $(\mathrm{size}(e)-1)$ ordinary edges
given by $\{v,e_*\}$ for $v$ in $e'$. Proposition \ref{propspine}
holds also for this construction except for the inequality of 
assertion (ii) which has to be replaced by the opposite inequality.
\end{rem}

%%%%%%%%%%%%%%%%%%%%%%%%%%%%%%%%%%%%%%%%%%%%%%%%%%%%%%%%%%%%%%%%%%%%%%%%%%%

\section{Glue maps, Pr\"ufer words and enumeration of 
labelled trees}\label{sectgluemap}

We consider a fixed rooted hypertree $T$ with root $r$ and non-root 
vertices $\mathcal V'$. 
We denote by $g:\mathcal V'\longrightarrow \mathcal V$ the map defined 
by $g(v)=w$ where $w$ is the unique neighbour of $v$ closer
to the root-vertex $r$ than $v$. The map $g$ sends thus every vertex 
$v$ of a reduced hyperedge $e'$ to the marked vertex $e_*$ of 
the associated hyperedge $e=e'\cup\{e_*\}$.
We extend $g$ to all vertices 
$\mathcal V=\mathcal V'\cup\{r\}$ of $T$ by setting $g(r)=r$.
We say that the map $g$ defines the marked vertices of $T$
and we call $g$ the \emph{glue-map} of the $r-$rooted tree $T$.

The sequence $v,g(v),g^2(v)=g(g(v)),\dots$ of iterates of $g$ is 
eventually constant and defines (up to repetitions of the root vertex)
the unique geodesic
joining a vertex $v$ of $T$ to the root vertex $r$.

Given a partition $\mathcal P$ of a set $\mathcal V'=\mathcal V
\setminus\{r\}$, a map $g:\mathcal V\longrightarrow \mathcal V$ 
is \emph{$\mathcal P-$admissible} if there exists an $r-$rooted
hypertree with vertices $\mathcal V$, Pr\"ufer
partition $\mathcal P$ and glue-map $g$.

\begin{prop} (i) Given a partition $\mathcal P$ of 
$\mathcal V'=\mathcal V\setminus\{r\}$,
a map $g:\mathcal V\longrightarrow\mathcal V$ is
$\mathcal P-$admissible if and only if the restriction of $g$ to 
a part of $\mathcal P$ is constant and for every vertex $v$ there exists 
an integer $k=k(v)$ such that
$g^k(v)=g^{k+1}(v)=r$ where $g^k$ denotes the 
$k-$fold iterate $g\circ g\circ \cdots\circ g$ of $g$.

\ \ (ii) $\mathcal P-$admissible partitions are in one-to-one correspondance
with $r-$rooted hypertrees having vertices $\mathcal V$ and Pr\"ufer
partition $\mathcal P$.
\end{prop} 

\noindent{\bf Proof} Associate to a part $e'$ of 
$\mathcal P$ the hyperedge with vertices $e'\cup
\{g(e')\}$. The condition $g^k(v)=r$ shows that these hyperedges
define a connected hypergraph. Since we have equality in 
inequality (\ref{edgedegineq}) for every finite connected subhypergraph
containing $r$, 
the resulting connected
hypergraph is a hypertree. This shows (i). 
Assertion (ii) is obvious.
\hfill$\Box$

\begin{rem} \label{remgluerestr}
A glue map $g$ of a hypertree $T$ with partition map $p$ is
completely determined by $g(r)=r$ and by its restriction 
$g\vert_{p(\mathcal V')}:p(\mathcal V')\longrightarrow\mathcal V$
to the set of fix-points (smallest elements in reduced hyperedges) of $p$.
\end{rem}

A \emph{Pr\"ufer word} is a one-to-one correspondence 
between the set of 
$\mathcal P-$admissible maps $\mathcal V'\longrightarrow\mathcal V$ 
and the set $\mathcal V^{k-1}$ of all words of length $k-1$ (with 
$k$ denoting the number of non-empty parts in the partition 
$\mathcal P$ of $\mathcal V'$) in the 
alphabet $\mathcal V$ satisfying the following two conditions:
\begin{enumerate}
\item{} The degree of a non-root vertex $v$ in the hypertree associated to a 
$\mathcal P-$admissible map $g$ is one more than the number of occurences
of the vertex $v$ in the word $W\in \mathcal V^{k-1}$ corresponding to $g$.
\item{} The Pr\"ufer word is given by a simple algorithm which is fast
(polynomial in any reasonable sense) for finite hypertrees.
\end{enumerate}

A Pr\"ufer word is thus an elegant way to recover the loss of information 
of the map $T\longmapsto \mathcal P(T)$ induced by the Pr\"ufer partition.

Formula (\ref{degreeequality}) implies that Condition (1) in a Pr\"ufer word
is also fulfilled by the root-vertex. Exempting the root from Condition (1) 
is motivated by Section \ref{sectinftrees} dealing with infinite hypergraphs.

A \emph{Pr\"ufer code} is a map $T\longmapsto (\mathcal P,W)$
where $\mathcal P$ is the Pr\"ufer partition of the non-root vertices 
$V'$ of a hypertree $T$ and where $W$ is a Pr\"ufer word.

The main result of this paper is a construction of a new Pr\"ufer 
word. This gives in particular a new proof of the following result:

\begin{thm} \label{thmpruf} Pr\"ufer codes exist.
\end{thm}

We give two proofs of Theorem \ref{thmpruf}. The first proof is
Selivanov's generalization of Pr\"ufer's construction, see \cite{Se}.
It amounts to the removal of hyperedges of a certain type
and it decreases the number of vertices. 

The second construction is the main result of this paper and seems
to be new. It is based on merging hyperedges with a common intersection
into one larger hyperedge and it does not change the set of vertices.
It has moreover the interesting feature that it works for suitably
defined infinite hypertrees, as outlined in Section
\ref{sectinftrees} and illustrated in Section \ref{sectexpleSn}.

We denote by $\mathcal T(\mathcal V,\mathcal P)$ the set of all 
finite rooted hypertrees with
vertices $\mathcal V=\mathcal V'\cup\{r\}$ and with a given 
fixed Pr\"ufer partition $\mathcal P$ of $\mathcal V'$.
Theorem \ref{thmpruf} implies easily the following standard 
result of enumerative combinatorics:

\begin{cor} Associating to a hypertree $T\in\mathcal T(\mathcal V,
\mathcal P)$ the monomial 
$$w(T)=\prod_{v\in\mathcal V}x_v^{\deg(v)-1}$$
we have
$$\sum_{T\in\mathcal T(\mathcal V,\mathcal P)}w(T)=\left(\sum_{v\in\mathcal V}
x_v\right)^{k-1}$$
where $k$ denotes the number of (non-empty) parts 
in $\mathcal P$.

In particular, there are 
$$S_2(k,n-1)n^{k-1}$$
hypertrees with $k$ hyperedges and vertices $\{1,\dots,n\}$
where $S_2(k,n)$ denotes the Stirling number of the second kind
enumerating the number of partitions of $\{1,\dots,n\}$
into $k$ non-empty subsets.
\end{cor}

%%%%%%%%%%%%%%%%%%%%%%%%%%%%%%%%%%%%%%%%%%%%%%%%%%%%%%%%%%%%%%%%%%%%%%%%%%%%%%%
\section{Hyperedges of leaf-type and the map 
$T\longmapsto W(T)$}\label{sectleaftype}

The next two sections deal only with finite hypertrees.

A hyperedge $e$ of a rooted hypertree $T$ is
of {\it leaf-type} if all vertices of the associated reduced 
hyperedge $e'$ are leaves.

\begin{prop} Every finite rooted hypertree not reduced to its root has a 
hyperedge of leaf-type.
\end{prop}

\noindent{\bf Proof} A hyperedge containing a vertex at maximal 
distance from the root-vertex is of leaf-type.\hfill$\Box$

\begin{prop} Given a hyperedge $e$ of a non-trivial 
hypertree $T$ with vertices $\mathcal V$, root $r$ and hyperedges $\mathcal E$, 
the set $\mathcal E\setminus\{e\}$ is the set of hyperedges of
an $r$-rooted hypertree with vertices $\mathcal V\setminus e'$
if and only if $e$ is of leaf-type.
\end{prop}

\noindent We leave the proof to the reader.\hfill$\Box$

We consider henceforth hypertrees with vertices given by a 
finite subset of $\mathbb  N$, rooted at the largest vertex and
with hyperedges totally ordered according to their smallest unmarked vertex.
We construct the \emph{Pr\"ufer word} $W(T)$ 
of such a hypertree $T$ by successively removing
the smallest hyperedge of leaf-type until reaching 
a trivial hypertree reduced to a unique hyperedge and by writing down the
sequence of marked vertices of the removed hyperedges. 

The following result is useful for the 
computation of the Pr\"ufer word of a tree given as a list
of hyperedges:

\begin{prop}\label{propsearchleaf}
A hyperedge $e$ of a rooted hypertree $T$ is of leaf-type
if and only if no element of the associated reduced 
hyperedge
$e'=e\setminus\{e_*\}$ occurs as a marked vertex 
among the other hyperedges of $T$.
\end{prop}

\noindent We leave the easy proof to the reader.\hfill$\Box$

\subsection{An example of a Pr\"ufer word}

The Pr\"ufer word $w_1\dots w_7$ of the hypertree $T$ represented by Figure 1 of Section \ref{sectexpleP} can be computed as follows: We start with 
the increasing sequence of all hyperedges, ordered according to their 
smallest non-marked vertex. We iterate then the following loop:
We search the first 
hyperedge $e$ of leaf-type using for example Proposition \ref{propsearchleaf}. 
We remove $e$ and we write down the marked vertex $e_*$ of the removed 
hyperedge $e$. We stop if only a unique hyperedge remains.

For our example represented in Figure 1, 
we get the increasing sequences of hyperedges
\begin{eqnarray*}
&&
\{1,10,12\}_8,{\it\{2\}_1},{\it\{3,9\}_4},\{4,7\}_8,
{\it\{5\}_{14}},{\it\{6\}_4},\{8,13\}_{14},{\it\{11\}_7}\\
&&{\it\{1,10,12\}_8},{\it\{3,9\}_4},\{4,7\}_8,
{\it\{5\}_{14}},{\it\{6\}_4},\{8,13\}_{14},{\it\{11\}_7}\\
&&{\it\{3,9\}_4},\{4,7\}_8,
{\it\{5\}_{14}},{\it\{6\}_4},\{8,13\}_{14},{\it\{11\}_7}\\
&&\{4,7\}_8,{\it\{5\}_{14}},{\it\{6\}_4},\{8,13\}_{14},{\it\{11\}_7}\\
&&\{4,7\}_8,{\it\{6\}_4},\{8,13\}_{14},{\it\{11\}_7}\\
&&\{4,7\}_8,\{8,13\}_{14},{\it\{11\}_7}\\
&&{\it\{4,7\}_8},\{8,13\}_{14}\\
&&{\it\{8,13\}_{14}}\\
\end{eqnarray*}
with hyperedges of leaf-type in italics.
The hypertree $T$ corresponds thus to the 
Pr\"ufer word $1\ 8\ 4\ 14\ 4\ 7\ 8$.
 %%%%%%%%%%%%%%%%%%%%%%%%%%%%%%%%%%%%%%%%%%%%%%%%%%%%%%%%%%%%%%%%%%%%%%%%%%%%%%%

\section{The inverse map $(\mathcal P,W)\longmapsto T$}
\label{sectinvW}

A part $e'$ in a partition $\mathcal P$ with $k$ non-empty parts 
of a subset $\mathcal S$ of $\{1,\dots,n-1\}$ 
is of \emph{leaf-type with respect to a word
$W\in \mathbb N^*$} over the alphabet $\mathbb N$ 
if $W$ involves no elements of $e'$. 

\begin{lem}\label{lempartleaftype} 
A partition $\mathcal P$ into  
$k$ non-empty parts of a subset of $\{1,\dots,n-1\}$ contains at 
least one part of leaf-type with respect to a word $W$ in the set
$\{1,\dots,n\}^{k-1}n$.
\end{lem}

{\bf Proof} The last letter $n$ of $W$ does not occur in any part of 
$\mathcal P$ and the number of remaining letters in $W$ is one 
less than the number of parts in $\mathcal P$.\hfill$\Box$

We consider a pair $(\mathcal P,W)$ consisting of a Pr\"ufer partition
$\mathcal P$
of $\mathcal S\subset \{1,\dots,n-1\}$ 
into $k$ non-empty parts and of a Pr\"ufer word 
$W\in(\mathcal S\cup\{n\})^{k-1}$. In order to construct the associated
hypertree $T$ rooted at $n$, it is enough to determine the glue map $g$ 
defining the marked vertex $e_*=g(e')$
of every reduced hyperedge $e'$ appearing in $\mathcal P$.
 This can be achieved as follows: We 
order the elements of $\mathcal P$ totally according to their smallest
element and we augment $W=w_1\dots w_{k-1}$ by
add a last letter $w_k=n$. We have thus $W=w_1\dots w_{k-1}n$. We iterate 
now the following loop: By Lemma \ref{lempartleaftype} there exists
a smallest part $e'$ of $\mathcal P$ which is of leaf-type 
with respect to $W=w_1\dots w_{k-1}$.
We get in this way the hyperedge $e'_{w_1}$ (given by all 
elements in $e'$ and by the marked vertex $w_1$) of $T$. 
We remove now $e'$ from $\mathcal P$, we erase $w_1$ in $W$ and we
iterate until $\mathcal P$ is empty.

\subsection{An example for the inverse map
$(\mathcal P,W)\longmapsto T$}

We reconstruct the hypertree $T$ of Figure 1 from its Pr\"ufer code
consisting of the Pr\"ufer partition
$$\{1,10,12\}\cup\{2\}\cup\{3,9\}\cup\{4,7\}\cup\{5\}\cup\{6\}\cup
\{8,13\}\cup\{11\}$$
and of the Pr\"ufer word $W=1\ 8\ 4\ 14\ 4\ 7\ 8$.
The computation for $T$ is as follows
$$\begin{array}{l|c}
\{1,10,12\},{\it\{2\}},{\it\{3,9\}},\{4,7\},{\it\{5\}},{\it\{6\}},
\{8,13\},{\it\{11\}}&1\\
{\it\{1,10,12\}},{\it\{3,9\}},\{4,7\},{\it\{5\}},{\it\{6\}},
\{8,13\},{\it\{11\}}&8\\
{\it\{3,9\}},\{4,7\},{\it\{5\}},{\it\{6\}},
\{8,13\},{\it\{11\}}&4\\
\{4,7\},{\it\{5\}},{\it\{6\}},
\{8,13\},{\it\{11\}}&14\\
\{4,7\},{\it\{6\}},
\{8,13\},{\it\{11\}}&4\\
\{4,7\},\{8,13\},{\it\{11\}}&7\\
{\it\{4,7\}},\{8,13\}&8\\
{\it\{8,13\}}&14\\
\end{array}$$
where the first columns displays relevant sets of reduced hyperedges with 
hyperedges of leaf-type in italics and where the last column
contains the letters of the Pr\"ufer word augmented with an additional 
letter $w_8=14$ representing the root-vertex. 
We get the hyperedges 
\begin{eqnarray*}
\{1,10,12\}_8,\{2\}_1,\{3,9\}_4,\{4,7\}_8,
\{5\}_{14},\{6\}_4,\{8,13\}_{14},\{11\}_7
\end{eqnarray*}
defining the hypertree $T$ of Figure 1 by
marking (indexing) the first reduced hyperedge 
of leaf type (written in italics) of every row with the corresponding letter of $W$.

%%%%%%%%%%%%%%%%%%%%%%%%%%%%%%%%%%%%%%%%%%%%%%%%%%%%%%%%%%%%%%%%%%%%%%%%%%%%%%%
\section{Star-reduction and the map
$T\longmapsto W^*(T)$}\label{sectstarred}

A \emph{hyperstar} is a hypertree of diameter at most $2$.
The \emph{center} of a hyperstar of diameter $2$ 
is the unique vertex adjacent to all other vertices.
It is given by the intersection of at least two hyperedges.
Every vertex is a center of the trivial hyperstar reduced to a unique
hyperedge.

The hyperstar $St(v)$ of a hypertree $T$ at a vertex $v$ of $T$ 
is the subtree of $T$ formed by $v$ and by all its neighbours.
Its hyperedges are all hyperedges of $T$ which contain $v$.

The {\it star-reduction} of $T$ at a vertex $v$ is the hypertree $*_v(T)$ 
obtained by replacing all hyperedges of $T$ involved in the hyperstar 
$St(v)$ by a unique hyperedge consisting of all vertices in $St(v)$.

\begin{prop} \label{lemstarcomm} (i) We have $*_v(T)=T$ if and only if
$v$ is leaf of $T$.

\ (ii) $v$ is a leaf of $*_v(T)$.

\ (iii) We have 
$*_v(*_w(T))=*_w(*_v(T))$ for any pair of vertices $v,w$ in $T$.
\end{prop}

Proofs are easy and left to the reader.

Assertion (iii) of Proposition \ref{lemstarcomm}
allows to define $*_S(T)$ for a subset $S$ of vertices.
Given a hypertree $T$ with vertices $\{1,\dots,n\}$,
we use the shorthand $*_{\leq v}(T)$ for the star-reduction 
$*_{\{1,\dots,v\}}(T)$ at 
the subset $\{1,\dots,v\}$ 
of all vertices not exceeding $v$.
All vertices $1,\dots,v$ of $*_{\leq v}(T)$ are leaves.
Similarly, we use $*_{<v}(T)$ for $*_{\{1,\dots,v-1\}}(T)$
using the convention $*_{<1}(T)=T$.

The Pr\"ufer word $W^*(T)$ of a hypertree $T$ with vertices $\{1,\dots,n\}$
and $k$ hyperedges is defined as follows: 

We set $W^*(T)=n^{k-1}$ if $T$ is a hyperstar with $k$ hyperedges
centered at its root-vertex $n$.

Otherwise, there exists a smallest non-leaf $v<n$ in $T$ and
we can define the increasing sequence (ordered by smallest
unmarked elements) $\mathcal E_v=(e_1,\dots,e_{k-1})$
of all $k-1$ hyperedges not containing the smallest non-leaf
$v$ as an unmarked vertex.
In other words, the sequence $\mathcal E_v$ is obtained by removing
the unique hyperedge $e$ containing $v$ in its reduced hyperedge $e'$
from the increasing sequence (ordered by smallest
unmarked elements) of all $k$ hyperedges of $T$. Since $v$
is a non-leaf there exist an increasing sequence 
$\mathcal S_v=1\leq i_1<\dots<i_a\leq k-1$
consisting of all $a=\deg(v)-1>0$ indices $i_1,\dots,i_a$
such that the hyperedges $e_{i_1},\dots,e_{i_a}$ of $\mathcal E_v$ 
have marked vertex $v$.
We set $w_{i_1}=\dots=w_{i_a}=v$ for the $a$ letters of $W^*(T)=w_1w_2\dots$
with indices $i_1,\dots,i_a$ in $\mathcal S_v$. 
The subword formed by the $k-1-a$ remaining letters 
$$w_1\dots w_{i_1-1}\widehat{w_{i_1}}w_{i_1+1}\dots \widehat{w_{i_2}}\dots
w_{i_a-1}\widehat{w_{i_a}}w_{i_a+1}\dots w_{k-1}$$
of $W^*(T)$ corresponds to the Pr\"ufer word $W^*(*_v(T))$
of the star reduction $*_v(T)$ of $T$ at the vertex $v$.

\begin{rem} The set $\mathcal T(n)$ of all hypertrees with vertices 
$1,\dots,n$ rooted at $n$ carries two interesting 
additional structures: 
\begin{enumerate}
\item{} It is a ranked poset for the order-relation given by 
$T\geq *_ST$ for any subset $S\subset \{1,\dots,n\}$ of vertices.
This poset has a unique minimal element given by the trivial hypertree
consisting of a unique hyperedge. Its rank function $rk(T)$ is given 
by the number $nl(T)$ of vertices which are non-leaves and its 
M\"obius function is $(-1)^{nl(T)}$.
\item{} The elements of $\mathcal T(n)$ are the vertices of a rooted hypertree.
The root vertex is again the trivial hypertree reduced to a unique hyperedge
with vertices $1,\dots,n$. The ancestor of a non-trivial 
tree $T$ is the star-reduction 
$s_a(T)$ with respect to the smallest vertex $a$ which is not a leaf.
\end{enumerate}
\end{rem}

%%%%%%%%%%%%%%%%%%%%%%%%%%%%%%%%%%%%%%%%%%%%%%%%%%%%%%%%%%%%%%%%%%%%%%%%%%%%%%%%

\subsection{An example for the
construction of the Pr\"ufer word $W^*(T)$}

We illustrate the computation of the Pr\"ufer word $W^*(T)$ by using 
once more our favourite tree with vertices $1,\dots,14$, root-vertex $14$ 
and hyperedges 
\begin{eqnarray*}
\{1,10,12\}_8,\{2\}_1,\{3,9\}_4,\{4,7\}_8,
\{5\}_{14},\{6\}_4,\{8,13\}_{14},\{11\}_7
\end{eqnarray*}
depicted in Figure 1 of Section \ref{sectexpleP}.

Non-leaves of $T$ in increasing order are $1,4,7,8,14$ and
increasing sequences of all hyperedges involved in non-trivial
star-reductions $*_{\leq v}(T)$ of $T$ are given by
$$\begin{array}{c|c}
v&\hbox{ hyperedges of }*_{<v}(T)\\
\hline
1&{\it\{1,10,12\}_8},\{2\}_1,\{3,9\}_4,\{4,7\}_8,
\{5\}_{14},\{6\}_4,\{8,13\}_{14},\{11\}_7,\\
4&\{1,2,10,12\}_8,\{3,9\}_4,{\it\{4,7\}_8},
\{5\}_{14},\{6\}_4,\{8,13\}_{14},\{11\}_7\\
7&\{1,2,10,12\}_8,{\it \{3,4,6,7,8,9\}_8},
\{5\}_{14},\{8,13\}_{14},\{11\}_7\\
8&\{1,2,10,12\}_8,\{3,4,6,7,9,11\}_8,
\{5\}_{14},{\it\{8,13\}_{14}}\\
&\{1,2,3,4,6,7,8,9,10,11,12,13\}_{14},\{5\}_{14}.\\
\end{array}
$$
The first column indicates the smallest non-leaf of $*_{<v}(T)$. 
The second column
consists of the increasing list of all hyperedges of $*_{<v}(T)$
with the hyperedge containing the smallest leaf as an unmarked 
vertex in italics.

Removing the italized hyperedges from the sequences in the second column,
we get the sequence $\mathcal E_v$ and the sequence 
$\mathcal S_v$ which determines the positions of the letters 
$1,4,7,8$ and $14$ in the Pr\"ufer word $W^*(T)=w_1\dots w_7$ of
$T$:
$$\begin{array}{c|c|c|c}
v&\mathcal E_v&\mathcal S_v&\mathcal E\setminus\mathcal E_v\\
\hline
1&\{2\}_1,\{3,9\}_4,\{4,7\}_8,
\{5\}_{14},\{6\}_4,\{8,13\}_{14},\{11\}_7&1&\{1,10,12\}_8\\
4&\{1,2,10,12\}_8,\{3,9\}_4,
\{5\}_{14},\{6\}_4,\{8,13\}_{14},\{11\}_7&2,4&\{4,7\}_8\\
7&\{1,2,10,12\}_8,\{5\}_{14},\{8,13\}_{14},\{11\}_7&4&\{3,\dots\}_8\\
8&\{1,2,10,12\}_8,\{3,4,6,7,9,11\}_8,
\{5\}_{14}&1,2&\{8,13\}_{14}\\
15&\{1,2,3,4,6,7,8,9,10,11,12,13\}_{14},\{5\}_{14}.&&\\
\end{array}$$
The necessary data for the final computation of $W^*(T)$ are summarized by
$$\begin{array}{c|c|cccccccccc}
v&\mathcal S_v&*_1&*_2&*_3&*_4&*_5&*_6&*_7\\
\hline
1&1&1&*_1&*_2&*_3&*_4&*_5&*_6\\
4&2,4&1&*_1&4&*_2&4&*_3&*_4\\
7&4&1&*_1&4&*_2&4&*_3&7\\
8&1,2&1&8&4&8&4&*_3&7\\
14&&1&8&4&8&4&14&7
\end{array}$$
and give $1\ 8\ 4\ 8\ 4\ 14\ 7$ for the Pr\"ufer word 
$W^*(T)$ of $T$.

From an algorithmic point of view it is perhaps more straightforward
to work with the partition map and with the glue map of $T$. 
Writing $p(w)_{g(w)}$ for the image of the vertex $w\in \{1,\dots,13\}$,
the partition map and the glue map 
of $*_{\leq v}(T)$ (with the convention $*_0(T)=T$) are given by
$$\begin{array}{c||c|c|c|c|c|c|c|c|c|c|c|c|c|}
v&1&2&3&4&5&6&7&8&9&10&11&12&13\\
\hline
0&1_8&2_1&3_4&4_8&5_{14}&6_4&4_8&8_{14}&3_4&1_8&11_7&1_8&8_{14}\\
1&1_8&1_8&   &   &     &   &  &     &   &1_8&&1_8&\\
4&   &   &3_8&3_8&&3_8&3_8&&3_8&&&&\\
7&&&&&&&&&&&3_8&&\\
8&1_{14}&1_{14}&1_{14}&1_{14}&&1_{14}&1_{14}&1_{14}&1_{14}&1_{14}&1_{14}&1_{14}&1_{14}\\
\end{array}$$
(unchanged values are omitted). 

Computing the partition map $p_i$ and the glue map $g_i$ of $*_{\leq i}(T)$
is straightforward by induction on $i$:
For $v$ a non-root vertex define $p_0(v)$ as the minimal element of the unique reduced hyperedge $e'$ containing $v$ and define $g_0(v)$ as the marked
vertex $e_*=e\setminus e'$ of the unique hyperedge $e$ 
whose associated reduced hyperedge $e'=E\setminus\{e_*\}\subset e$ 
contains $v$. 

We suppose now $p_{i-1}$ and $g_{i-1}$ constructed.
We consider
$$a=\min\left(p_{i-1}(i),\min_{v\in\mathcal V',g_{i-1}(v)=i}p_{i-1}(v)\right)$$
and we set $p_i(v)=a$ if either $p_{i-1}(v)=p_{i-1}(i)$ or
$g_{i-1}(v)=i$. We leave $p_i(v)=p_{i-1}(v)$ unchanged otherwise, ie.
if $p_{i-1}(v)\not=p_{i-1}(i)$ and $g_{i-1}(v)\not=i$.
The integer $a$ is of course the minimal element in the unique
reduced hyperedge of $*_{\leq i}(T)$ which contains $i$.

We set $g_i(v)=g_{i-1}(v)$ if $g_{i-1}(v)\not=i$ and 
we set $g_i(v)=g_{i-1}(i)$ if $g_{i-1}(v)=i$. Otherwise stated,
the marked vertex $g_{i-1}(i)$ of the unique 
reduced hyperedge $e'$ in $*_{<i}(T)$ which contains
$i$ is not affected by star-reduction at $i$ except if it is equal 
to $i$. In this case it is replaced by the marked vertex $g_{i-1}(i)$ 
of the unique reduced hyperedge in $*_{<i}(T)$ which contains $i$.

Using Remark \ref{remgluerestr} we can condense
the computations for $\mathcal S_v$ to
$$\begin{array}{c|c|c|c}
v&&\mathcal S_v\\
\hline
1&2_1,3_4,4_8,5_{14},6_4,8_{14},11_7&1&1_8\\
4&1_8,3_4,5_{14},6_4,8_{14},11_7&2,4&4_8\\
7&1_8,5_{14},11_7&3&3_8\\
8&1_8,3_8,5_{14}&1,2&8_{14}\\
\end{array}$$
by choosing smallest unmarked representatives in hyperedges.

%%%%%%%%%%%%%%%%%%%%%%%%%%%%%%%%%%%%%%%%%%%%%%%%%%%%%%%%%%%%%%%%%%%%%%%%%%%%%%

\section{The inverse map $(\mathcal P,W^*)\longmapsto 
T$}\label{sectinvW*}

Given a Pr\"ufer code $(\mathcal P,W^*)$ where $\mathcal P$
is a partition of $\{1,\dots,n-1\}$ into $k$ non-empty parts and 
where $W^*\in\{1,\dots,n\}^{k-1}$ is a word of length $k-1$ with 
letters in $\{1,\dots,n\}$, there exists a unique hypertree 
$T$ such that $\mathcal P=\mathcal P(T)$ is the Pr\"ufer partition 
of $T$ and $W^*=W^*(T)$, defined by the construction of Section
\ref{sectstarred}, is the Pr\"ufer word of $T$. 

If $W^*=n^{k-1}$, the pair $(\mathcal P,W^*)$ is Pr\"ufer code of
the hyperstar centered at the root-vertex $n$ with
reduced hyperedges given by the parts of $\mathcal P$.
Otherwise there exists a smallest letter $v<n$ occuring with 
strictly positive multiplicity
$a>0$ in $W^*$. We denote by $\mathcal P_v=(s_1,\dots,s_{k-1})$ 
the increasing sequence (ordered with respect to smallest elements)
of all $k-1$ parts not containing the vertex $v$ of $\mathcal P$. 
If $i_1,\dots,i_a$ are the $a$ indices of all letters equal to $v$
in the word $W^*=w_1\dots w_{k-1}$ then the parts $s_{i_1},\dots,s_{i_a}$
correspond to all reduced hyperedges of a (not yet constructed) 
hypertree $T$ with marked vertex $v$.
Denoting by $e'_v$ the unique part of $\mathcal P$ containing $v$,
we consider the partition $\tilde{\mathcal P}$
obtained by merging the $a+1$ parts $s_{i_1},\dots,s_{i_a}$
and $e'_v$ into a larger part $\tilde{e'}$. 
We denote by $\widetilde{W^*}$ the word of length
$k-1-a$ obtained by removing all $a$ letters equal to $v$ from $W^*$.
The Pr\"ufer word $\widetilde{W^*}$ of the pair
$(\widetilde{\mathcal P},\widetilde{W^*})$ contains no letter $\leq a$.
It is thus by descending induction on $a$ the Pr\"ufer code of 
a unique hypertree $\widetilde T$ 
with $\{1,\dots,a\}$ contained in the set of leaves. 
More precisely, the recursively defined hypertree $\widetilde T$ 
is the star-reduction $*_v(T)$ at $v$  
of the hypertree $T$ corresponding to $(\mathcal P,W^*)$.
The hyperedge $\tilde e$ associated to the part $\tilde{e}'$
of $\widetilde T$ splits into $a+1$ hyperedges of $T$ in the obvious way:
$a$ hyperedges with marked vertex $v$ have reduced hyperedges 
$s_{i_1},\dots,s_{i_a}$. The marked vertex of the hyperedge corresponding to
the last part $e'_v$ involved in $\tilde{e'}$ is given by the marked vertex 
of the hyperedge associated to $\tilde{e'}$ in $\tilde T=*_v(T)$. This defines
the hypertree $T$ uniquely.

%%%%%%%%%%%%%%%%%%%%%%%%%%%%%%%%%%%%%%%%%%%%%%%%%%%%%%%%%%%%%%%%%%%%%%%%%%

\subsection{An example for the inverse map}
We reconstruct the hypertree $T$ of Figure 1
from its Pr\"ufer code $(\mathcal P,W^*)$ consisting of the 
Pr\"ufer partition 
\eject
$$
\mathcal P=\{1,10,12\}_8,\{2\}_1,\{3,9\}_4,\{4,7\}_8,
\{5\}_{14},\{6\}_4,\{8,13\}_{14},\{11\}_7$$
(with parts totally ordered by minimal elements)
and of the Pr\"ufer word $W^*=1\ 8\ 4\ 8\ 4\ 14\ 7$.

We have 
$$\begin{array}{c|l}
1&{\it\{1,10,12\}},\{2\}_1,\{3,9\},\{4,13\},\{5\},\{6\}_4,\{8,13\},
\{11\}\\
4&\{1,2,10,12\},\{3,9\}_4,{\it\{4,7\}},\{5\},\{6\}_4,\{8,13\},\{11\}\\
7&\{1,2,10,12\},{\it\{3,4,6,7,9\}},\{5\},\{8,13\},\{11\}_7\\
8&\{1,2,10,12\}_8,\{3,4,6,7,9,11\}_8,\{5\},{\it\{8,13\}}\\
15&\{1,2,3,4,6,7,8,9,10,11,12,13\}_{14},\{5\}_{14}
\end{array}$$
The indices $1,4,7,8$ are added according to the positions
of the letters $1,4,7,8$ and in the words
\begin{eqnarray*}
W^*&=&1\ 8\ 4\ 8\ 4\ 14\ 7\\
W^*\setminus\{2\}&=&8\ 4\ 8\ 4\ 14\ 7\\
W^*\setminus\{2,4\}&=&8\ 8\ 14\ 7\\
W^*\setminus\{2,4,7\}&=&8\ 8\ 14\\
\end{eqnarray*}
after removal of the italicized
part containing the index under consideration. Reduced hyperedges
of the last row are all marked by the root $14$. 

Parts in every row are completely
ordered according to smallest elements and are obtained
from the parts of the previous row by merging all parts involving the 
vertex considered in the previous row (and by copying the remaining parts).

The marked vertex of a reduced hyperedge $e'$ is now given by the index of 
the first indexed superset $\tilde e'\supset e'$ encountered when moving down 
the rows. We get thus the hyperedges
$$\{1,10,12\}_8,\{2\}_1,\{3,9\}_4,\{4,7\}_8,\{5\}_{14},
\{6\}_4,\{8,13\}_{13},\{11\}_7$$
of our favourite hypertree $T$ depicted in Figure 1.

The following table illustrates the algorithm 
by giving partition maps for $*_{\leq v}(T)$ 
and by giving partial glue-maps (denoted by 
$p(v)_{g(v)}$, see Remark \ref{remgluerestr})
$$\begin{array}{|c||c|c|c|c|c|c|c|c|c|c|c|c|c|c|}
\hline
v&1&2&3&4&5&6&7&8&9&10&11&12&13&\\
\hline
p_0&1&2&3&4&5&6&4&8&3&1&11&1&8&\\
\hline
1&\multicolumn{13}{c|}{2_1,3,4,5,6,8,11}&1\\
\hline
p_1&1&1&3&4&5&6&4&8&3&1&11&1&8&\\
\hline
4&\multicolumn{13}{c|}{1,3_4,5,6_4,8,11}&4\\
\hline
p_4&1&1&3&3&5&3&3&8&3&1&11&1&8&\\
\hline
7&\multicolumn{13}{c|}{1,5,8,11_7}&3\\
\hline
p_7&1&1&3&3&5&3&3&8&3&1&3&1&8&\\
\hline
8&\multicolumn{13}{c|}{1_8,3_8,5}&8\\
\hline
p_8&1&1&1&1&5&1&1&1&1&1&1&1&1&\\
\hline
&\multicolumn{13}{c|}{1_{14},5_{14}}&\\
\hline
\end{array}
$$
The table encodes the information for the glue map $g$ as follows:
Suppose we want to determine the marked vertex $g(7)$ of the reduced hyperedge
containing vertex $7$. The second row (denoted by $p_0$) contains 
the information for the partition map of $*_{\leq 0}(T)=T$. It
shows that the smallest element in the reduced hyperedge (of $T$) 
containing $7$ is
$3$. This implies $g(7)=g(4)$ and we are reduced to compute $g(4)$.
Nothing interesting happens to vertex $4$ during the star-reduction at 
vertex $1$. After that, vertex $4$ is involved in the star-reduction at vertex 
$4$ and becomes an element of the reduced hyperedge of $*_{\leq 4}(T)$ 
with smallest element $p_4(4)=3$. We switch thus our attention to the vertex 
$3$. The next row indicates that the hyperedge containing $3$ of
$*_{\leq 4}(T)$ has marked vertex $8$. We have thus $g(7)=8$
for the value $g(7)$ of the glue map $g$ at the vertex $8$.

Proceeding similarly we get the complete information 
$$\begin{array}{|c||c|c|c|c|c|c|c|c|c|c|c|c|c|}
\hline
v&1&2&3&4&5&6&7&8&9&10&11&12&13\\
\hline
g&8&1&4&8&14&4&4&14&4&8&7&8&14\\
\hline
\end{array}$$
for the glue map $g$ of the tree $T$ depicted in Figure 1.

%%%%%%%%%%%%%%%%%%%%%%%%%%%%%%%%%%%%%%%%%%%%%%%%%%%%%%%%%%%%%%%%%%%%%%%%%%%%%%

\section{Infinite hypertrees}\label{sectinftrees}

The construction of the Pr\"ufer code $(\mathcal P,W^*)$ 
based on mergings of hyperedges
works perfectly well for infinite hypertrees with vertices $\{1,2,3,\dots\}
\cup\{\infty\}$ rooted at $\infty$. It encodes such a hypertree $T$ with 
infinitely many hyperedges by a Pr\"ufer partition of $\{1,2,\dots\}$
into infinitely
many non-empty parts and an infinite Pr\"ufer word $W^*=w_1w_2\dots\in 
(\mathbb N\cup\{\infty\})^{\mathbb N}$ with an arbitrary vertex $v$ (which can be the root 
vertex) of $T$ occuring $\deg(v)-1$ times in $W^*$ where $\deg(v)$ can be 
infinite.

The Pr\"ufer map is however not onto: A pair $(\mathcal P,W^*)$ consisting 
of a partition of $\{1,2,\dots\}$ into
infinitely many parts and an infinite word $W^*\in 
(\mathbb N\cup\{\infty\})^{\mathbb N}$ corresponds in general to no infinite 
rooted hypertree. In order to have a one-to-one correspondence, we introduce
in this section ideally rooted 
hypertrees with vertices $\{1,2,\dots\}\cup\{\infty\}$.
Such objects are hyperforests having at most one component which is
an ordinary (finite or infinite) rooted hypertree together with an 
arbitrary large (and perhaps infinite) number of infinite trees with marked
ends playing the role of the root vertex $\infty$.

Observe that infinite hypertrees with vertices $\{1,2,\dots\}$ rooted
at $\infty$ which have only finitely many hyperedges are essentially 
the same as finite hypertrees from the point of view of the Pr\"ufer word
$W^*$. We leave the easy discussion for this class of rooted hypertrees to the 
reader.

\subsection{Ends of hypertrees and ideally rooted hypertrees}

Two infinite geodesics $\gamma,\gamma':\mathbb N\longrightarrow \mathcal V$ 
of an infinite hypergraph $G$ are \emph{equivalent} if 
$d(\gamma(n),\gamma'(n))$ is ultimately constant. Equivalence classes
of such infinite geodesics are called \emph{ends} of $G$.

An \emph{ideally rooted hypertree} is a hyperforest with a choice of an
end in every connected component not containing the root vertex.
We call the connected component containing the root of an ideally
rooted hypertree the \emph{root component}. The root component can be finite
(and perhaps reduced to its root) or infinite. All other components are
\emph{ideal components}. They contain always infinitely many hyperedges.

An ideally rooted hypertree has a marked vertex $e_*$ in every hyperedge.
The marked vertex $e_*$ of a hyperedge $e$ in the root component is defined in
the usual way as the unique vertex of $e$ which is closest to the root.
The marked vertex $e_*$ of a hyperedge $e$ in an ideal component $C$ is defined
as the unique vertex closest to $\gamma(n)$ for $n$ huge enough where 
$\gamma:\mathbb N\longrightarrow \mathcal V$ is a fixed geodesic defining 
the equivalence class of the marked end of $C$. We leave it to the reader 
to show that $e_*$ is well defined and depends only on the equivalence class
of $\gamma$.

\subsection{Partition maps and glue maps of ideally rooted hypertrees}

Partition maps of ideally rooted hypertrees are idempotent lowering maps of 
the set $\mathbb N*=\mathbb N\setminus \{0\}$ into itself. Glue maps
are maps of the set $\mathbb N*\cup\{\infty\}$  
admitting the fixpoint $g(\infty)=\infty$ as their only recurrent element.
Extending partition maps by $p(\infty)=\infty$, 
pairs maps $p,g:\mathbb N^*\cup \{\infty\}\longrightarrow
\mathbb N^*\cup \{\infty\}$ fixing $\infty=p(\infty)=g(\infty)$
formed by a lowering idempotent map $p$ and a map $g$ with
$\infty=g(\infty)$ as its unique recurrent element correspond
to partition maps and the glue maps of ideally rooted trees if and 
only if $g=g\circ p$.

\subsection{The Pr\"ufer code of an ideally rooted hypertree}

Proofs are straightforward and omitted in this informal Section.

The Pr\"ufer partition $\mathcal P=\mathcal P(T)$ of an ideally
rooted hypertree $T$ is defined in the obvious way as the partition of 
the set $\mathcal V\setminus\{r\}$ of non-root vertices
with parts $e\setminus\{e_*\}$ given by all reduced hyperedges.

The \emph{glue-map} of an ideally rooted hypertree $T$ is the map $g:
\mathcal V\longrightarrow \mathcal V$ having the root $r=g(r)$ as its unique
fixpoint and given by $g(v)=e_*$ for a non-root vertex $v$ arising as an 
unmarked element of the hyperedge $e$.

The Pr\"ufer word $W^*(T)$ of an ideally rooted hypertree $T$ with vertices
$\{1,2,\dots\}\cup\{\infty\}$ rooted at $\infty$ is well-defined and
given by an infinite word $w_1w_2\dots$ with a finite letter 
$n\in \mathbb N$ occuring exactly $\deg(n)-1$ times. The degree
$\deg(n)$ of a vertex $n$ can be finite or infinite.
The letter $\infty$ corresponding to the root vertex occurs at most 
$\deg(\infty)-1$ times in $W^*(T)$ where $\deg(\infty)\in \{0,1,2,\dots\}
\cup\{\infty\}$ is defined as 
the degree of the vertex $\infty$ in the root-component.

The exact number of occurences of $\infty$ in $W^*(T)$ can be strictly 
smaller than $\deg(\infty)-1$. More precisely, we order the connected 
components of $T\setminus\{\infty\}$ according to their smallest vertex.
Denoting by $\widetilde{C}$ the smallest ideal component, we erase 
all connected components $\geq \widetilde{C}$ from $T$ and we denote by 
$\widetilde{\deg}(\infty)$ the degree of $\infty$ in the resulting 
ordinary rooted 
hypertree. The number of occurences of the letter $\infty$ in $W^*(T)$
is then given by $\max(0,\widetilde{\deg}(\infty)-1)$. 

We have now a one-to-one correspondence between the set of 
ideally rooted hypertrees
having vertices $\{1,2,\dots\}\cup\{\infty\}$ rooted at $\infty$ and having
infinitely many hyperedges and the set of Pr\"ufer codes consisting of 
a Pr\"ufer partition of $\mathbb N$ into infinitely many non-empty parts
and an infinite Pr\"ufer word which can be an arbitrary element of 
$\left\{\{1,2,\dots,\}\cup\{\infty\}\right\}^{\mathbb N}$.

Ideally rooted hypertrees (with infinitely many hyperedges)
which are locally finite correspond to Pr\"ufer codes with Pr\"ufer
partitions involving only finite parts and with Pr\"ufer words 
involving all letters $\left\{\{1,2,\dots,\}\cup\{\infty\}\right\}$
with finite multiplicity.

\section{An example giving rise to bijections of $S_n$}
\label{sectexpleSn}

The simplest infinite tree with vertices $\mathbb N\cup\{\infty\}$,
rooted at $\infty$, is given by a halfline originating at the root 
$\infty$ with vertices $v_1,v_2,\dots \in \mathbb N$ at distance $1,2,\dots$
of the root-vertex $\infty$. Setting $v_0=\infty$, each vertex $v_i$ other
than the root-vertex $v_0$ has thus exactly two neighbours 
$v_{i-1}$ and $v_{i+1}$.
The root vertex $\infty$ has a unique neighbour $v_1$. Such a tree is 
completely described by the permutation $i\longmapsto v_i$ of the set 
$\{1,2,3,\dots\}$ and every permutation $\sigma$ of $\{1,2,\dots\}$
describes a unique such tree. The Pr\"ufer word
$W^*$ of such a tree yields a permutation $\psi$ of $\{1,2,\dots\}$.
(Caution: not every permutation of $\{1,2,\dots\}$ corresponds to such a 
tree: most permutations give rise to trees with ideal components.)
The Pr\"ufer partition of such a tree is of course the trivial partition
of $\{1,2,\dots\}$ into singletons and thus carries no information.

A particularly nice subset of permutations is given by so-called
\lq\lq finitely-supported'' permutations
moving only finitely many elements of the infinite
set $\{1,2,\dots\}$. Such a permutation $\sigma$ satisfies $\sigma(m)=m$
for every integer $m$ larger than some natural integer $n$ and thus can
be considered as an element of the finite permutation group
$S_n$ acting in the usual way on $\{1,\dots,n\}$. It is easy to see 
that the Pr\"ufer word $W^*$ of such a tree has this property again.
Thus the Pr\"ufer word defines a bijection of $S_n$ which respects
the the obvious inclusion of $S_{n-1}$ in $S_n$ as the subset of 
all permutations fixing $n$.

We describe now this map for $n\leq 4$. We write 
$\left(\begin{array}{cccc}\sigma(1)&\sigma(2)&\dots&\sigma(n)\end{array}
\right)$ for a permutation $i\longmapsto \sigma(i)$ of $\{1,2,\dots,n\}$.

In the case $n=1$ there is a unique permutation. It fixes every element
of $\{1,2,\dots\}$ and the associated Pr\"ufer word $W^*$ is again
the identity permutation.

The unique non-trivial permutation $\sigma$ in $S_2$ (extendend to 
a permutation of $\{1,2,3,\dots\}$ by setting $\sigma(i)=i$ for all $i>2$)
gives again rise to $W^*(\sigma)=\sigma$.

In the case $n=3$, the image $W^*(\sigma)$ of $\sigma$ is already known
for the two permutations of $S_2\subset S_3$. The remaining four
permutations form two orbits defined by the image $\sigma(3)$ of the 
largest integer $3$. Note that we have always $\psi(i)=\sigma(i)=i$
for the partition $\psi$ encoded by the Pr\"ufer word $\psi=W^*(\sigma)$
of a partition $\sigma$ such that $\sigma(i)=i$ for all $i>n$.

In the case $n=4$, the map $\sigma\longmapsto W^*(\sigma)$ gives 
rise to two orbits
$$
\left(\begin{array}{cccc}
2&3&4&1\end{array}\right),
\left(\begin{array}{cccc}
4&2&3&1\end{array}\right),
\left(\begin{array}{cccc}
3&2&4&1\end{array}\right),
\left(\begin{array}{cccc}
4&3&2&1\end{array}\right)$$
and 
$$
\left(\begin{array}{cccc}
2&4&3&1\end{array}\right),
\left(\begin{array}{cccc}
3&4&2&1\end{array}\right)$$
associated to permutations such that $\sigma(4)=1$. We have
finally one orbit
\begin{eqnarray*}
&\left(\begin{array}{cccc}
1&3&4&2\end{array}\right),
\left(\begin{array}{cccc}
4&1&3&2\end{array}\right)
\left(\begin{array}{cccc}
3&1&4&2\end{array}\right),&\\
&\left(\begin{array}{cccc}
3&4&1&2\end{array}\right),
\left(\begin{array}{cccc}
1&4&3&2\end{array}\right),
\left(\begin{array}{cccc}
4&3&1&2\end{array}\right)&
\end{eqnarray*}
associated to all permutations such that $\sigma(4)=2$ and one
orbit
\begin{eqnarray*}
&&\left(\begin{array}{cccc}
1&2&4&3\end{array}\right),
\left(\begin{array}{cccc}
1&4&2&3\end{array}\right),
\left(\begin{array}{cccc}
4&2&1&3\end{array}\right),\\
&&
\left(\begin{array}{cccc}
2&1&4&3\end{array}\right),
\left(\begin{array}{cccc}
2&4&1&3\end{array}\right),
\left(\begin{array}{cccc}
4&1&2&3\end{array}\right)
\end{eqnarray*}
consisting of all permutations with $\sigma(4)=3$.

\noindent Roland BACHER, Universit\'e Grenoble I, CNRS UMR 5582, Institut 
Fourier, 100 rue des maths, BP 74, F-38402 St. Martin d'H\`eres, France.

\noindent e-mail: Roland.Bacher@ujf-grenoble.fr

\end{document}